\documentclass{amsart}
\usepackage{amssymb}
\let\p\pi
\let\q\rho
\def\calA{\mathcal{A}}
\def\calB{\mathcal{B}}
\def\calC{\mathcal{C}}
\def\P{\mathcal{P}}
\def\vecbeta{\vec\beta}
\def\betaelt_#1^#2{\beta_{#1}^{(#2)}}
\let\realbigwedge\bigwedge
\def\bigwedge{\realbigwedge\nolimits} 
\def\bwedge{{\textstyle\bigwedge}}
\newtheorem{theorem}{Theorem}[section]
\newtheorem{defn}[theorem]{Definition}
\newtheorem{lemma}[theorem]{Lemma}

\def\C{k}
\def\k{p}
\def\Z{{\mathbb Z}}
\def\Spec{\mathop{\rm Spec}\nolimits}

\def\Char{\mathop{\rm char}\nolimits}
\let\ve\varepsilon
\def\Ker{\mathop{\rm Ker}\nolimits}

\title{Universality of
Rank 6  Pl\"ucker Relations and Grassmann Cone Preserving Maps}

\author{Alex Kasman}
\author{Kathryn Pedings}
\author{Amy Reiszl}
\author{Takahiro Shiota}

\begin{document}

\begin{abstract}  The Pl\"ucker relations define a projective embedding
  of the Grassmann variety $Gr(\k,n)$.
  We give another finite set of quadratic equations
  which defines the same embedding, and whose elements all have rank 6.
  This is achieved by constructing a certain finite set of linear maps
  $\bigwedge^\k\C^n\to\bigwedge^2\C^4$,
  and pulling back the unique Pl\"ucker relation on $\bigwedge^2\C^4$.
  We also give a quadratic equation depending on $(\k+2)$ parameters
  having the same properties.
\end{abstract}

\maketitle
% \markright{\uppercase{Rank 6 Pl\"ucker Relations and GCP Maps}}


\section{Introduction}
\subsection{Pl\"ucker Relations}
Throughout, let $\C$ be a field, and let
$e_1$, \dots, $e_n$ be a basis of the vector space $\C^n$. Define the
coordinates $\{\Pi_{i_1,\dots,i_\k}\}_{1\le i_1<\cdots<i_\k\le n}$ on
$\bigwedge^\k\C^n$ by
\begin{equation}
\bigwedge^\k\C^n\ni
\omega=\sum_{1\leq i_1<i_2<\cdots<i_\k\leq n} \Pi_{i_1i_2\ldots i_\k}
e_{i_1}\wedge e_{i_2}\wedge\cdots\wedge e_{i_\k},\label{omega}
\end{equation}
and extend them to arbitrary indices in $\{1,\dots,n\}^\k$ by making them
antisymmetric.
% $\Pi_{i_1\ldots i_ji_{j+1}\ldots i_\k}=-\Pi_{i_1\ldots i_{j+1}i_j\ldots i_\k}$,
% and $\Pi_{i_1\ldots i_ji_j\ldots i_\k}=0$.

An element $\omega\in\bigwedge^\k\C^n$ is called \textit{decomposable}
if it can be written in the form
$\omega=v_1\wedge v_2\wedge\cdots\wedge v_\k$ for some $v_i\in\C^n$;
otherwise it is called \textit{indecomposable}.
The \textit{Grassmann cone} $\Gamma^\k\C^n=\{\omega\in\bigwedge^\k\C^n\mid
\omega=v_1\wedge v_2\wedge\cdots\wedge v_\k\hbox{ for some }v_i\in\C^n\}$
is the set of decomposable elements in $\bigwedge^\k\C^n$.
The \textit{Pl\"ucker relations}
\cite{Bourbaki,GriffithsHarris,HodgePedoe,KL}
\begin{equation}
P_{A,B}(\omega):=\sum_{i=1}^{\k+1}(-1)^{i-1}
\Pi_{a_1a_2\dots a_{\k-1}b_i}\Pi_{b_1b_2\dots b_{\k+1}\setminus b_i}=0,
\label{DefPlucker}
\end{equation}
where
$A=\{a_1,\dots,a_{\k-1}\}$,~$B=\{b_1,\dots,b_{\k+1}\}\subset\{1,\dots,n\}$,
and where the $\setminus b_i$ at the end of the indices indicates the absence of
$b_i$ from the indices, hold if and only if $\omega\in\Gamma^\k\C^n$,
making $\Gamma^\k\C^n$ a $\C^\times$-invariant affine variety.\footnote{
 To be precise, while \eqref{DefPlucker} cut out $\Gamma^\k\C^n$
 set-theoretically for any $\C$ \cite[Theorem 1]{KL},
 more general Pl\"ucker relations may be needed to define it
 scheme-theoretically if $\Char(\C)>0$ \cite{Abe}.
 In this paper we shall concentrate on the set-theoretic aspect, so
 the simpler Pl\"ucker relations as in \eqref{DefPlucker} suffice.
}
The quotient
$(\Gamma^\k\C^n\setminus\{0\})/\C^\times\subset\mathbb P(\bigwedge^\k\C^n)$
is the Grassmann variety.

We do not need to consider all the choices of indices
$A$ and $B$: Since rearranging the elements of $A$ or $B$ only affects
$P_{A,B}$ by total change in sign, it suffices to consider
$A$ and $B$ whose elements are listed in increasing order.
Moreover, if $A\subset B$ then $P_{A,B}=0$, and
if $A\setminus(A\cap B)$ is a one point set $\{a\}$ then exchanging $a$
with any element of $B\setminus(A\cap B)$ only affects
$P_{A,B}$ by total change in sign.
So we take
$$
  \P(\k,n)=\left\{P_{A,B}\left|\
  \begin{array}{@{}l@{}}
  A,B\subset \{1,\ldots,n\},\ A=\{a_1,\ldots,a_{\k-1}\},\
  B=\{b_1,\ldots,b_{\k+1}\}\\
  \hbox{with $a_1<\cdots<a_{\k-1}$ and $b_1<\cdots<b_{\k+1}$},\ A\not\subset B\hbox{, and}\\
  \hbox{if $A\setminus(A\cap B)=\{a\}$ then
  $a<b$ for any $b\in B\setminus(A\cap B)$}
  \end{array}
  \right.\right\}
$$
as a set of generators of Pl\"ucker relations, and it suffices to define
the Grassmann cone:
$\Gamma^\k\C^n=\{\omega\in\bigwedge^\k\C^n\mid P(\omega)=0\hbox{ for all }
P\in\P(\k,n)\}$.

By definition we have $\P(\k,n)=\emptyset$ if $\min\{\k,n-\k\}\le1$.
The first nontrivial case $(\k,n)=(2,4)$ yields
$\P(2,4)=\{P_{\{1\}\{2,3,4\}}\}$, so
$\omega=\sum_{1\leq i<j\leq 4} \Pi_{ij} e_i\wedge e_j\in\bigwedge^2\C^4$
is decomposable if and only if
\begin{equation}
P_{\{1\}\{2,3,4\}}(\omega)=
\Pi_{12}\Pi_{34}-\Pi_{13}\Pi_{24}+\Pi_{14}\Pi_{23}=0.\label{gr24}
\end{equation}
% For example $\P(2,4)=\{P_{\{1\}\{2,3,4\}}\}$, so
% $\omega=\sum_{1\leq i<j\leq 4} \Pi_{ij} e_i\wedge e_j\in\bigwedge^2\C^4$
% is decomposable if and only if
% \begin{equation}
% P_{\{1\}\{2,3,4\}}(\omega)=
% \Pi_{12}\Pi_{34}-\Pi_{13}\Pi_{24}+\Pi_{14}\Pi_{23}=0.\label{gr24}
% \end{equation}
% This is the simplest nontrivial
% Pl\"ucker relation % \cite[Chapter 6]{GriffithsHarris}
% since $\P(\k,n)=\emptyset$ if $\min\{\k,n-\k\}\le1$.

% \subsection{The Quadric Rank of the Pl\"ucker Relations}

The \textit{rank} of a quadratic form is the rank of the symmetric matrix
which defines it (unless $\Char(\C)=2$, in which case
a different definition of rank is used \cite{Chevalley}).
The rank of $P_{A,B}$ as a quadratic form on $\bigwedge^\k\C^n$
is twice the number $|B\setminus(A\cap B)|$ of nonvanishing terms
in \eqref{DefPlucker}.
So the set $\P(\k,n)$ consists of quadratic forms of every even rank
from 6 up to $2\min\{\k,n-\k\}+2$, and the
Pl\"ucker relations in $\P(\k,n)$ all have rank 6 only when $\min\{\k,n-\k\}=2$.
The literature on algebraic geometry occasionally demonstrates an interest in the rank of the Pl\"ucker relations \cite{rk6ref1,rk6ref2}, with particular attention paid to the simplest, rank 6 case.

\subsection{Grassmann Cone Preserving Maps}\label{sec:GCP}

A linear map $G\colon\bigwedge^\k\C^n\to\bigwedge^{\k'}\C^{n'}$
is said to be a Grassmann cone preserving map (GCP map for short) if
$$
G(\Gamma^\k\C^n)\subset \Gamma^{\k'}\C^{n'}.
$$

The \textit{induced map}
$\bwedge^\k L \colon\bigwedge^\k\C^n\to\bigwedge^\k\C^{n'}$ of a linear map
$L\colon\C^n\to\C^{n'}$ is a GCP map.
Another GCP map is the \textit{dual isomorphism}
\begin{eqnarray*}
\delta\colon\bigwedge^\k\C^n&\to&\bigwedge^{n-\k}\C^n\\
e_{i_1}\wedge\cdots\wedge e_{i_\k}&\mapsto& e_{j_1}\wedge\cdots\wedge e_{j_{n-\k}}
\end{eqnarray*}
where $i_1<\cdots<i_\k$, $j_1<\cdots<j_{n-\k}$ and
$\{i_1,\ldots,i_\k\}\cup \{j_1,\ldots,j_{n-\k}\}=\{1,\ldots,n\}$.
% $j_m\neq i_l$ for $1\leq l\leq \k$, $1\leq m\leq n-\k$.

% Aside from trivial GCP maps (those whose image contains only decomposable elements),
Every nontrivial GCP map (those which do {\it not\/} send the whole of $\bigwedge^\k\C^n$
to decomposables) can be written as a composition of maps of these
two types \cite{Westwick3}.

\subsection{Motivation}

In the theory of classical integrable systems, it is well-known
that the KP hierarchy of soliton equations, written in an appropriate (Hirota) form,
is nothing but the Pl\"ucker relations for an infinite dimensional Grassmannian in the space of functions \cite{Sato}.
Somewhat remarkably,
% in this case there is a 3-term (i.e., rank 6) quadratic relation with parameters which serves as a generating function for a large enough set of equations to define the same Grassmannian \cite{Duzhin,GK,Miwa,Tak-Tak}.
in that setting a single 3-term (i.e., rank 6) %% HERE
quadratic functional equation with parameters suffices to encode the  
entire hierarchy \cite{GK,Miwa,Tak-Tak}.
However, the literature on algebraic geometry does not appear to have a similar result on the ``universality" of the 3-term Pl\"ucker relation \eqref{gr24}, i.e., that the 3-term relation is in a sense the only Pl\"ucker relation one needs.

In this paper, we will show such universality
% purely in the framework of multilinear algebra,
by pulling back \eqref{gr24} by various GCP maps
to obtain a finite set of rank 6 equations on
$\bigwedge^\k\C^n$ which suffices to cut out the cone of decomposables
$\Gamma^\k\C^n\subset\bigwedge^\k\C^n$ set-theoretically.
% as $\P(\k,n)$ does.
%Thus the indecomposability of $\omega\in\bigwedge^\k\C^n$ is always
%visible when it is mapped in a suitable way to $\bigwedge^2\C^4$,
%the simplest case.
Inspired by the soliton theory observation, we will also construct
a parameter-dependent rank 6 equation which determines
the decomposability of $\omega$.

\section{A Set of Polynomials with Quadric Rank 6}

We will define another set of quadratic forms $\P'(\k,n)$ to be used later to
cut out the Grassmann cone $\Gamma^\k\C^n$ just as $\P(\k,n)$ does.
First we introduce a convenient notation:
%which naturally reflects the multilinearity of wedge products:
for a 2-vector $\vec i_j=(i_j^{(0)},i_j^{(1)})\in\{1,\dots,n\}^2$, let
$$
\Pi_{i_1i_2 \ldots i_{j-1}\vec i_j i_{j+1}\ldots i_\k}=
\Pi_{i_1i_2 \ldots i_{j-1}\,i_{j}^{(0)}\,i_{j+1}\ldots i_\k}+
\Pi_{i_1i_2 \ldots i_{j-1}\,i_{j}^{(1)}\,i_{j+1}\ldots i_\k},
$$
and extend it inductively to the case where two or more indices are 2-vectors.

\begin{defn}\label{pprime-def}
For $\calA=\{\alpha_1,\alpha_2,\alpha_3,\alpha_4\}\subset\{1,\dots,n\}$,
$\calB=\{\vecbeta_1,\vecbeta_2,\ldots,\vecbeta_m\}$ with
$\vecbeta_{i}=(\betaelt_{i}^{0},\betaelt_{i}^{1})\in\{1,\dots,n\}^2$, and
$\calC=\{\gamma_1,\ldots,\gamma_{\k-m-2}\}\subset\{1,\dots,n\}$, let
\begin{equation}
P'_{\calA,\calB,\calC}=\pi_{12}\pi_{34}-\pi_{13}\pi_{24}+\pi_{14}\pi_{23},
\quad\hbox{where}\quad
\pi_{ij}=
\Pi_{\alpha_i\alpha_j\vecbeta_1\ldots\vecbeta_m\gamma_{1}\ldots\gamma_{\k-2-m}}.
\label{PprimeABC}
\end{equation}
\end{defn}

Note that $P'_{\calA,\calB,\calC}$ is nothing but 
\eqref{gr24} with the six variables $\Pi_{ij}$ replaced by $\pi_{ij}$.  
%(One may think of $P'_{\calA,\calB,\calC}$ is as a rank 6 element of $\P(\k,n)$ combined with a change of basis of
%the underlying vector space $\C^n$, as we will see below.)
% In fact, $P'_{\calA,\calB,\calC}$ is the same as $P_{A,B}$ with $A=\{\betaelt_1^1,\betaelt_2^1,\ldots,\betaelt_m^1,\gamma_1,\ldots,\gamma_{\k-2-m},\alpha_1\}$ and $B=\{\alpha_2,\alpha_3,\alpha_4,\betaelt_1^1,\betaelt_2^1,\ldots,\betaelt_m^1,\gamma_1,\ldots,\gamma_{\k-2-m}\}$ combined with the change of basis that adds $e_{\betaelt_i^0}$ to $e_{\betaelt_i^1}$ for each $i=1,\ldots,m$.
Thus, it has rank at most six
as a quadratic form of the original coordinates $\Pi_{i_1\dots i_p}$. Moreover,
if $\alpha_i$, $\beta_j^{(\mu)}$ and $\gamma_l$ are all distinct it has rank exactly six.

\begin{lemma}\label{pprime-expand}Each $P'_{\calA,\calB,\calC}$ in
\eqref{PprimeABC} %, not necessarily subject to condition \eqref{OrderABC},
can be written as a $\Z$-linear combination of elements of
%$\{\pm P\mid P\in\P(\k,n)\}$:
$\P(\k,n)$:
$$
P'_{\calA,\calB,\calC}=\sum_{\vec\mu,\vec\nu\in\{0,1\}^m} P_{\{\betaelt_1^{\mu_1}\ldots\betaelt_m^{\mu_m}\gamma_1\ldots \gamma_{\k-m-2}\alpha_1\}\{\alpha_2\alpha_3\alpha_4\betaelt_1^{\nu_1}\ldots\betaelt_m^{\nu_m}\gamma_1\ldots \gamma_{\k-m-2}\}}.
$$
\end{lemma}
\begin{proof}
For each fixed choice of $\vec\mu$ and $\vec\nu$ the first three terms in expansion \eqref{DefPlucker} of $P$ are exactly those that appear in \eqref{PprimeABC}.
The next $m$ terms, i.e., the terms involving $\betaelt_i^{\nu_i}$,
are zero if $\mu_i=\nu_i$, and cancel if $\mu_i\ne\nu_i$ with the term where $\mu_i$ and $\nu_i$ are switched.
The last $\k-m-2$ terms, the terms involving $\gamma_i$, are all zero.
\end{proof}

The lemma implies that if some $P'_{\calA,\calB,\calC}(\omega)\neq0$,
then $\omega$ is indecomposable.
Conversely, in Section~\ref{sec:pp2} we will see that some 
$P'_{\calA,\calB,\calC}(\omega)\ne0$ 
if $\omega$ is indecomposable.
However, 
\textit{not} every element of $\P(\k,n)$ is a
linear combination of rank 6 quadratic forms $P'_{\calA,\calB,\calC}$. %\footnote{
For instance, a
tedious but straightforward computation shows that the Pl\"ucker relation $P_{\{1,2,3\},\{4,5,6,7,8\}}$ cannot be  a linear combination of  $P'_{\calA,\calB,\calC}$. % }

The quadratic forms $P'_{\calA,\calB,\calC}$ are not independent.  Clearly, rearranging the elements of each index set affects them by at most a change in sign.
% Those with the same index appearing more than once in $\calA$, $\calB$ or $\calC$ can be achieved as linear combinations
If the indices $\alpha_i$, $\betaelt_j^\ve$ and $\gamma_l$ are not all distinct,
then $P'_{\calA,\calB,\calC}$ is a linear combination
of those with distinct indices.  There are other,
less obvious linear relations like
$$
\sum_{i\in\Z/3\Z}P'_{\{\nu_i,\alpha_1,\alpha_2,\alpha_3\},
\calB\cup\{(\nu_{i+1},\nu_{i+2})\},\calC}
-\!\!\sum_{0\le i,j\le2;i\ne j}P'_{\{\nu_i,\alpha_1,\alpha_2,\alpha_3\},
\calB,\calC\cup\{\nu_j\}}=0.
$$
However, instead of using these identities to find a basis for the linear span
of the set of all $P'_{\calA,\calB,\calC}$, let us introduce a subset just
suited for our purpose of set-theoretic characterization of the Grassmann cone.

%Consequently, in analogy with our definition of  $\P(\k,n)$ above, we will identify a relatively small\footnote{}  subset, $\P'(\k,n)$, of these polynomials whose indices satisfy certain conditions but will be suitable for our purposes.

\begin{defn}\label{pprime-set-def}
Let $\P'(\k,n)$ be the set of all $P'_{\calA,\calB,\calC}$ where the triples
$\calA, \calB, \calC$ satisfy
\begin{equation}
\left\{
\begin{array}{@{\,}l@{}}
\smash{\hbox{the indices $\alpha_i$, $\beta_j^{(\mu)}$ and $\gamma_l$
are all distinct and such that}}\\[2.5pt]
\alpha_i<\alpha_{i+1},\ \
\betaelt_i^0<\betaelt_i^1,\ \
\betaelt_i^\nu<\betaelt_{i+1}^\nu,\ \
\betaelt_i^0<\alpha_1,\ \
\betaelt_i^1<\alpha_3,\ \
\gamma_i<\gamma_{i+1}.
\end{array}
\right.
\label{OrderABC}
\end{equation}
\end{defn}

It is a simple exercise in combinatorics to find the number of elements
of $\P'(\k,n)$:
$$
|\P'(\k,n)|=\sum_{m=0}^M
\frac{n!}{(2m+4)!(\k-m-2)!(n-\k-m-2)!}(C_{m+2}-C_{m+1}),
$$
where $M=\min\{\k,n-\k\}-2$, and
$C_r=\frac1{r+1}\binom{2r}r$ is the Catalan number \cite{Stanley}.
Rewriting this, we see that in general
$\P'(\k,n)$ is a much smaller set than $\P(\k,n)$:
$$
|\P'(\k,n)|=\sum_{m=0}^M\frac{3m+3}{(2m+4)(2m+3)}a_m
\le
\frac14a_0+\sum_{m=1}^Ma_m=|\P(\k,n)|,
$$
where the equality holds if and only if $M\le0$, and where
$$
a_m:=\frac{n!}{(m+1)!(m+3)!(\k-m-2)!(n-\k-m-2)!}
$$
is the number of pairs of subsets $A$,~$B\subset\{1,\dots,n\}$ such that
$|A|=\k-1$, $|B|=\k+1$ and $|A\cap B|=\k-m-2$, which is also
the number of the rank $2m+6$ elements in $\P(\k,n)$ if $m>0$,
and is four times this number if $m=0$.

\section{Decomposability and GCP Maps to $\Gamma^2\C^4$}

% If $\omega\in\bigwedge^\k\C^n$ is decomposable then for every GCP map $G\colon\bigwedge^\k\C^n\to\bigwedge^2\C^4$ it is true that $G(\omega)$ is decomposable as well.
In this section we will construct a finite set of GCP maps from $\bigwedge^\k\C^n$ to $\bigwedge^2\C^4$,
indexed by the same triples $\calA$, $\calB$, $\calC$ as in
Definition~\ref{pprime-set-def}, such that
if $\omega$ is indecomposable then for \textit{some} GCP map $G$ in this set $G(\omega)$ is also indecomposable.
First we define a linear map from $\C^n$ to $\C^{\k+2}$ determined by these indices and prove a lemma addressing a question of vector geometry.

For any $S\subset\{1,\dots,n\}$ let
$\C^S=\bigoplus_{i\in S}\C e_i\subset \C^n$.
For any $S$,~$T\subset\{1,\dots,n\}$ and a map $f\colon S\to T$, define
$f^*\colon k^T\to k^S$ and $f_*\colon k^S\to k^T$ by
$\sum_{j\in T}a_je_j\mapsto\sum_{i\in S}a_{f(i)}e_i$ and
$\sum_{i\in S}a_ie_i\mapsto\sum_{i\in S}a_ie_{f(i)}$, respectively.
Thus if $i_S$ is the natural inclusion $S\subset\{1,\dots,n\}$,
then $\pi_S:=i_S^*\colon\C^n\to\C^S$ is the projection
$\sum_{i=1}^n a_ie_i\mapsto\sum_{i\in S}a_ie_i$.
Writing $S=\{\xi_1,\dots,\xi_{|S|}\}$ with $\xi_1<\cdots<\xi_{|S|}$,
we also define
$\tau_S\colon S\to\{1,\dots,|S|\}$ by $\tau_S(\xi_i)=i$, so that
$\tau_{S*}\colon\C^S\stackrel\sim\to\C^{|S|}$, $\tau_{S*}(e_{\xi_i})=e_i$,
is an isomorphism.

\begin{defn}\label{X-def}
Suppose $\calA$, $\calB$ and $\calC$ are as in Definition~\ref{pprime-def},
subject to condition \eqref{OrderABC}.
Let $B_\ve=\{\betaelt_1^\ve,\dots,\betaelt_m^\ve\}$ ($\ve=0$,~$1$), and let
$S'=\calA\cup B_1\cup\calC$ and $S=S'\cup B_0$.
Let $\varphi\colon B_0\stackrel\sim\to B_1$ be the order-preserving map
$\varphi(\betaelt_i^0)=\betaelt_i^1$, $i=1,\dots,m$, and extend it to a map
from $S$ to $S'$, still denoted by $\varphi$, by $\varphi(i)=i$ outside $B_0$.
Noting that $|S'|=\k+2$, we let
% $\{\betaelt_{1}^{1},\ldots,\betaelt_{m}^{1}\}$ written in increasing order so that $\xi_i<\xi_{i+1}$.
$X:=X_{\calA,\calB,\calC}\colon\C^n\to\C^{\k+2}$ be the composition of linear
maps
$\tau_{S'*}\circ\varphi_*\circ\pi_S\colon\C^n\to\C^S\to\C^{S'}\to\C^{\k+2}$,
% $\C^n\stackrel{\pi_S}\to\C^S\stackrel{\varphi_*}\to\C^{S'}\stackrel{\tau_{S'*}}\to\C^{\k+2}$,
i.e., writing $S'=\{\xi_1,\dots,\xi_{|S'|}\}$ with $\xi_1<\cdots<\xi_{|S'|}$
we have
$$
X(e_i)=\left\{\begin{array}{@{\,}ll}
e_j& \smash{\hbox{if $i=\xi_j$,
or if $\exists l\ [i=\betaelt_l^0$ and $\betaelt_l^1=\xi_j]$}},\\[3pt]
% X(e_{\betaelt_{j}^{1}})& \hbox{if $i=\betaelt_{j}^{0}$},\cr
0&\hbox{otherwise}.
\end{array}\right.
$$
\end{defn}

\begin{lemma}\label{lem:subspaces}
Let $V_0$ and $V_1$ be $\k$-dimensional subspaces of $\C^n$.  If
$q:=q(V_0,V_1):=\k-\dim(V_0\cap V_1)=\dim(V_0+V_1)-\k$ satisfies $q\geq 2$,
there exist $\calA$, $\calB$ and $\calC$ subject to \eqref{OrderABC} so that 
\begin{itemize}
\item $\hat V_0=X_{\calA,\calB,\calC}(V_0)$ and $\hat V_1=X_{\calA,\calB,\calC}(V_1)$ are $\k$-dimensional subspaces of $\C^{\k+2}$,
\item  $\dim (\hat V_0\cap \hat V_1)=\k-2$,
\item $\{X_{\calA,\calB,\calC}(e_{\alpha_i})\mid i=1,2,3,4\}$ is linearly independent modulo $(\hat V_0\cap \hat V_1)$.
\end{itemize}
\end{lemma}

\begin{proof}
Take a minimal $S\subset\{1,\dots,n\}$
such that $\pi_S|_{V_\ve}\colon V_\ve\to \C^S$ ($\ve=0$,~1) are injective, and
such that $q_S:=q(\pi_S(V_0),\pi_S(V_1))\ge2$.
Let $V_\ve':=\pi_S(V_\ve)$.
We have
$$
q_S\le q
$$
since $V_0'\cap V_1'\supset\pi_S(V_0\cap V_1)\simeq V_0\cap V_1$
implies $\k-q_S\ge \k-q$,
and
\begin{equation}
\k+q_S=|S|
\label{foo1}
\end{equation}
since the minimality of $S$ implies $V_0'+V_1'=\pi_S(V_0+V_1)=\C^S$.

If $q_S=2$ so that $\dim \C^S/(V_0'\cap V_1')=4$, choose
$\calA\subset S$ such that $|\calA|=4$ and the composition of the
maps $\C^\calA\hookrightarrow \C^S\to \C^S/(V_0'\cap V_1')$
is an isomorphism,
and let $\calB=\emptyset$ and $\calC=S\setminus\calA$.
It is easy to see that $X_{\calA,\calB,\calC}$
satisfies the desired properties.

If $q_S>2$, then for any $i\in S$ the minimality of $S$ implies that
$\pi_{S\setminus\{i\}}|_{V_0}$ or $\pi_{S\setminus\{i\}}|_{V_1}$
cannot be injective (otherwise we would have $q_{S\setminus\{i\}}=q_S-1\ge2$), thus
\begin{equation}
e_i\in V_0'\quad\hbox{or}\quad e_i\in V_1'\quad(\forall i\in S).
\label{foo2}
\end{equation}
Take any $B_0$, $B_1\subset S$ such that
\begin{equation}
\C^{B_\ve}\subset V_\ve'\quad\hbox{and}\quad
\C^{B_\ve}\cap V_{1-\ve}'=\{0\}\quad(\ve=0,1).
\label{foo3}
\end{equation}
Note that \eqref{foo3} implies $B_0\cap B_1=\emptyset$ and
$|B_\ve|\le q_S$.
Let $\delta=0$ if $|B_0|\le|B_1|$; otherwise let $\delta=1$.
Let $S'=S\setminus B_\delta$.  Take any inclusion
$\varphi\colon B_\delta\hookrightarrow B_{1-\delta}\subset S'$, and
extend it to a map from $S$ to $S'$
by letting $\varphi(i)=i$ for $i\in S'$.
As a consequence  of \eqref{foo3} we have $\Ker\varphi_*\cap V_\ve'=\{0\}$, so that
$V_\ve'':=\varphi_*(V_\ve')$ are $\k$-dimensional subspaces of $\C^{S'}$.
Since $V_0''+V_1''=\varphi_*(V_0'+V_1')=\varphi_*(\C^S)=\C^{S'}$
has dimension $|S|-|B_\delta|=\k+q_S-|B_\delta|$ (see \eqref{foo1}), we
have
\begin{equation}
q(V_0'',V_1'')=q_S-|B_\delta|=q_S-\min\{|B_0|,|B_1|\}.
\label{count_q}
\end{equation}

If we further assume that $(B_0,B_1)$ is a maximal pair satisfying \eqref{foo3},
then $|B_0|=|B_1|=q_S$.  To prove this,
it suffices to show that $B_\delta$ is not maximal if $|B_\delta|<q_S$.
By \eqref{foo2} we have
$$
e_i\in V_0''\quad\hbox{or}\quad e_i\in V_1''\quad(\forall i\in S'),
$$
but since $|\{i\in S'\mid e_i\in V_\ve''\}|\le\dim(V_\ve'')=\k$
($\ve=0$, 1) and $|S'|=\k+q_S-|B_\delta|$,
if $|B_\delta|<q_S$ then there exists $i\in S'$ such that
\begin{equation}
e_i\not\in V_{1-\delta}''.
\label{foo5}
\end{equation}
Let $B_\delta'=B_\delta\cup\{i\}$.
Since $i\in S'=S\setminus B_\delta$ we have $i\not\in B_\delta$,
so that $B_\delta'\supsetneq B_\delta$.
By \eqref{foo5} we have $e_i\not\in V_{1-\delta}'$; hence
by \eqref{foo2} we have $e_i\in V_\delta'$.
This and the first formula in \eqref{foo3} (with $\ve=\delta$) yield
\begin{equation}
\C^{B_\delta'}=\C^{B_\delta}\oplus \C e_i\subset V_\delta'.
\label{foo6}
\end{equation}
Using \eqref{foo5}, the first formula in \eqref{foo3} with $\ve=1-\delta$
and the second formula in \eqref{foo3} with $\ve=\delta$, we have
$$
e_i\not\in\varphi_*^{-1}(V_{1-\delta}'')=
V_{1-\delta}'\oplus\bigoplus_{i\in B_\delta}\C(e_i-e_{\varphi(i)})
=V_{1-\delta}'\oplus \C^{B_\delta},
$$
hence
$
\C^{B_\delta'}\cap V_{1-\delta}'=\{0\},
$
which together with \eqref{foo6} shows that $B_\delta$ is
not maximal.

Take maximal $B_\ve$, denote them by $\overline B_\ve$, and write
$\overline B_{\ve}=\{b_1^{\ve},\ldots,b_{q_S}^{\ve}\}$ with
$b_1^{\ve}<\cdots<b_{q_S}^{\ve}$. We then let
$B_\ve=\{b_1^{\ve},\ldots,b_{q_S-2}^{\ve}\}\subset\overline B_\ve$,
\begin{gather*}
\calA=\{b_{q_S-1}^{0},b_{q_S}^{0},b_{q_S-1}^{1},b_{q_S}^{1}\}
=(\overline B_0\setminus B_0)\cup(\overline B_1\setminus B_1),
\\
\calB=\{(b_i^{0},b_i^{1})\mid i=1,\dots,q_S-2\}
\quad\hbox{and}\quad
\calC=S\setminus(\overline B_0\cup\overline B_1)=S\setminus(\calA\cup B_0\cup B_1).
\end{gather*}
As subsets of $\overline B_\ve$, $B_\ve$ satisfy \eqref{foo3}.
Taking order-preserving
$\varphi\colon B_0\stackrel\sim\to B_1$ in the above argument
(note $\delta=0$ since $|B_0|=|B_1|$)
we obtain $\varphi\colon S\to S'$ and
$X:=X_{\calA,\calB,\calC}=\tau_{S'*}\circ\varphi_*\circ\pi_S$
as in Definition~\ref{X-def}. The map $X$
is regular on each $V_{\ve}$ (so that $\hat V_\ve=\tau_{S'*}V_\ve''$ has dimension $\k$),
and $X(e_{\alpha_i})=\tau_{S'*}(e_{\alpha_i})$ remain linearly independent modulo $\hat V_0\cap\hat V_1=\tau_{S'*}(V_0''\cap V_1'')$.
% Furthermore, the identification of $q_S-2$ basis elements in $V_0$ only with corresponding ones in $V_1$ only increases the dimension of the intersection so that
Finally, applying \eqref{count_q} to $B_\ve$ we have
$q(V_0'',V_1'')=q_S-|B_\delta|=2$, and
$\dim(\hat V_0\cap\hat V_1)=\dim(V_0''\cap V_1'')=\k-2$ as desired.
\end{proof}

Now, let us define a GCP map associated to the same indices:
\begin{defn}\label{G-def}
The Grassmann cone preserving map $G_{\calA,\calB,\calC}\colon\bigwedge^\k\C^n\to\bigwedge^2\C^4$ is defined as
$$
G_{\calA,\calB,\calC}=\bwedge^2 Z\circ \delta\circ\bwedge^\k X_{\calA,\calB,\calC},
$$
where
$X_{\calA,\calB,\calC}$ is as in Definition~\ref{X-def},
$\delta\colon\bigwedge^\k\C^{\k+2}\to\bigwedge^2\C^{\k+2}$ is the dual isomorphism, and $Z\colon\C^{\k+2}\to\C^4$ is the linear map
$$Z(e_i)=\left\{\begin{array}{@{\,}ll}
e_j & \smash{\hbox{if $e_i=X_{\calA,\calB,\calC}(e_{\alpha_j})$},}\\[2pt]
0& \hbox{otherwise}.
\end{array}\right.
$$
\end{defn}

\begin{theorem}\label{thm:existsABC}
The element $\omega\in\bigwedge^\k\C^n$ is indecomposable if and only if there exists a choice of indices $\calA$, $\calB$, $\calC$ so that $G_{\calA,\calB,\calC}(\omega)\in\bigwedge^2\C^4$ is indecomposable.  
\end{theorem}
\begin{proof}
Since $G_{\calA,\calB,\calC}$ is a GCP map, if $\omega$ is decomposable
$G_{\calA,\calB,\calC}(\omega)$ is also decomposable for any choice of
$\calA$, $\calB$, and $\calC$.
So, let us assume that $\omega$ is indecomposable and show that for an appropriate choice of the indices, its image in $\bigwedge^2\C^4$ is also indecomposable.

Suppose first that $\omega$ can be written as a sum of \textit{two} decomposable elements
\begin{equation}
\omega=v_1\wedge\cdots\wedge v_\k + w_1\wedge\cdots\wedge w_\k. \label{SumOfTwo}
\end{equation}
Then $\omega$ is indecomposable if and only if \cite{Westwick3} the $\k$-dimensional subspaces $V_0=\langle v_1,\ldots,v_\k\rangle$ and $V_1=\langle w_1,\ldots,w_\k\rangle$ of $\C^n$ satisfy $\dim(V_0\cap V_1)\le \k-2$.  Applying Lemma~\ref{lem:subspaces} gives us a choice of $\calA$,
$\calB$, $\calC$ such that 
$$
(\bwedge^\k X_{\calA,\calB,\calC})(\omega)=\omega_1 \wedge \omega_2\in \bigwedge^\k\C^{\k+2},
$$
where $\omega_1$ is an indecomposable element of the second exterior power of the 4-dimensional space
$\C^{A'}:=X_{\calA,\calB,\calC}(\C^\calA)$, and
% $0\ne\omega_2\in\bigwedge^{\k-2}(\hat V_0\cap\hat V_1)$.
$\omega_2$ is a nonzero element of $\bigwedge^{\k-2}(\hat V_0\cap\hat V_1)$.
Since $\C^{A'}\cap(\hat V_0\cap\hat V_1)=\{0\}$, letting %%HERE
$A'':=\{1,\dots,\k+2\}\setminus A'$ we have
$\omega_{2,0}:=(\bwedge^{\k-2}\pi_{A''})\omega_2\ne0$,
which implies
$\bwedge^2\pi_{A'}\circ\delta(\pi_1\wedge\pi_2)=
\bwedge^2\pi_{A'}\circ\delta(\pi_1\wedge\pi_{2,0})=c\delta_{A'}(\pi_1)$,
where $c\ne0$ and $\delta_{A'}$ is the dual isomorphism on $\bigwedge^*\C^{A'}$.
Since $Z=\tau_{A'*}\circ\pi_{A'}$, this implies
% =e_{j_1}\wedge\cdots\wedge e_{j_{\k-2}}$ is simply the wedge product of the $\k-2$ basis elements other than the four indexed by the $\alpha_i$'s.
$G_{\calA,\calB,\calC}(\omega)=\bwedge^2Z\circ\delta(\omega_1\wedge\omega_2)
=(\bwedge^2\tau_{A'*})(c\delta_{A'}(\omega_1))\ne0$. 

This proves the claim when $\omega$ can be written as a sum of two decomposable elements.  This is always true in $\bigwedge^2\C^4$, but not in general.  We proceed by induction on the dimension $n$ with the case $n=4$ as our initialization.

% It uniquely determines $\omega_1\in\bigwedge^{\k-1}\C^n$ and $\omega_2\in\bigwedge^\k\C^n$ to require that 

Regarding $\C^{n-1}$ as the subspace $\langle e_1,\ldots,e_{n-1}\rangle$ of
$\C^n=\langle e_1,\ldots,e_n\rangle$, let $\omega_1\in\bigwedge^{\k-1}\C^{n-1}$
and $\omega_2\in\bigwedge^\k\C^{n-1}$ be such that
$
\omega=\omega_1\wedge e_n + \omega_2.
$
%  with the condition that $\omega_i\wedge e_n=0$ only if $\omega_i=0$.
Now, if $\omega_1$ is indecomposable, we make use of the induction hypothesis on $\bigwedge^{\k-1}\C^{n-1}$ to get $\calA$, $\calB$, and $\calC$ so that $G_{\calA,\calB,\calC}(\omega_1)$ is indecomposable.  If we consider instead $\calC'=\calC\cup\{n\}$ which has cardinality one greater than $\calC$ then $G_{\calA,\calB,\calC'}(\omega_1\wedge e_n)=G_{\calA,\calB,\calC}(\omega_1)$ and $G_{\calA,\calB,\calC'}(\omega_2)=0$ since $\bwedge X(\omega_2)=0$.  Consequently, $\calA,\calB,\calC'$ satisfies the claim.

On the other hand, if $\omega_2$ is indecomposable we make use of the induction hypothesis to obtain $\calA$, $\calB$ and $\calC$.  In this case, without any further modification we have that $G(\omega)=G(\omega_2)$ is indecomposable  since the absence of the number $n$ in the indices will result in $\bwedge X$ which annihilates $\omega_1\wedge e_n$.

The only other possibility is that both $\omega_1$ and $\omega_2$ are decomposable, which returns us to the case that was proved initially.
\end{proof}

\section{Determining Decomposability Using the Elements of $\P'(\k,n)$}\label{sec:pp2}

Applying Theorem~\ref{thm:existsABC},
we now show that the rank 6 quadratic forms in $\P'(\k,n)$ -- like the Pl\"ucker relations in $\P(\k,n)$ -- are capable of characterizing the decomposables.

\begin{theorem}
For any $\omega\in\bigwedge^\k\C^n$, all the elements of $\P'(\k,n)$ vanish at $\omega$ if and only if $\omega$ is decomposable.
\end{theorem}
\begin{proof}
By Theorem~\ref{thm:existsABC}, $\omega$ is decomposable if and only if
$G_{\calA,\calB,\calC}(\omega)\in\bigwedge^2\C^4$ is decomposable for all choices of  $\calA$, $\calB$, $\calC$ as in Definition~\ref{pprime-set-def} and $G_{\calA,\calB,\calC}$ as in Definition~\ref{G-def}.
However, the decomposability of an element of $\bigwedge^2\C^4$ is determined by the single Pl\"ucker relation \eqref{gr24}.

So it is sufficient to note that substituting the 6 coordinates of $G_{\calA,\calB,\calC}(\omega)$ into \eqref{gr24} yields precisely the equation $P'_{\calA,\calB,\calC}=0$.  To see this, note that writing
$$
G_{\calA,\calB,\calC}(\omega)=\sum_{1\leq i_1<i_2\leq 4}\hat\Pi_{i_1i_2}e_{i_1}\wedge e_{i_2}
$$
we can compute $\hat\Pi_{i_1i_2}$ directly.  It is necessarily independent of $\Pi_{j_1j_2j_3\ldots j_\k}$ unless exactly two of $\alpha_1,\ldots,\alpha_4$, one each of $\betaelt_j^0$ and $\betaelt_j^1$ for each $1\leq j\leq m$ and all of $\gamma_1,\ldots,\gamma_{\k-m-2}$ are represented amongst $j_1,\ldots,j_\k$ since the wedge products of which these are the coefficients are in the kernel of $G$.  Moreover, of the remaining elements the only restriction on those which arise in $\hat \Pi_{i_1i_2}$ is that $\alpha_{i_1}$ and $\alpha_{i_2}$ are \textit{not}
present amongst the indices.  In particular,
$$
\hat\Pi_{i_1i_2}=(\epsilon_{ij})\Pi_{\alpha_1\alpha_2\alpha_3\alpha_4\vecbeta_1\ldots\vecbeta_m\gamma_{1}\ldots\gamma_{\k-2-m}\backslash\alpha_{i_1}\alpha_{i_2}}
$$
where $\epsilon_{ij}\in\{-1,1\}$  and substituting these into \eqref{gr24} results in $\pm P'_{\calA,\calB,\calC}$.
\end{proof}

\section{Parameter Dependent Formulation}
\label{sect:ParamDep}

In this section we will introduce a GCP map depending polynomially on the
$\k+2$ free parameters $x_1,\dots,x_{\k+2}$, and use it to
determine the decomposability of $\omega\in\bigwedge^\k\C^n$.
Let $R=\C[x_1,\dots,x_{\k+2}]$, and let $Q$ be the field of fractions of $R$.
We call $\eta\in\bigwedge^a R^b$ decomposable if
$\eta\otimes_RQ\in\bigwedge^aQ^b$ is decomposable, i.e., if
$\eta\otimes_RQ\in\Gamma^aQ^b$.
If $\eta$ is decomposable, then for any $\mathfrak p\in\Spec R$
there exist $v_1,\dots,v_{a}\in(R_{\mathfrak p})^{b}$ such that
$\eta=v_1\wedge\cdots\wedge v_{a}$.
Here $R_{\mathfrak p}$ is the localization of $R$ at $\mathfrak p$.

\begin{defn} Let $X$ be the $(\k+2)\times n$ matrix whose
$(i,j)$ entry is $x_i^{j-1}$,
regarded as a $Q$-linear map $X\colon Q^n\to Q^{\k+2}$.
Let $Z=(I\ O)$ be the $4\times(\k+2)$ matrix consisting
of $4\times4$ unit matrix and $4\times(\k-2)$ zero matrix,
also regarded as a $Q$-linear map $Z\colon Q^{\k+2}\to Q^4$.
The GCP map $\bwedge^2 Z\circ\delta\circ\bwedge^\k X\colon
\bigwedge^\k Q^n\to\bigwedge^2Q^4$ restricted to
$\bigwedge^\k\C^n\subset \bigwedge^\k Q^n$
defines $G\colon\bigwedge^\k\C^n\to\bigwedge^2R^4$.  Let
$H=G^*P_{\{1\}\{2,3,4\}}\colon\bigwedge^\k\C^n\to R$, i.e.,
for $\omega\in\bigwedge^\k\C^n$, writing
$$
G(\omega)=\sum_{1\leq i<j\leq 4}\hat\Pi_{ij}e_i\wedge e_j,\quad
\hat\Pi_{ij}\in R,
$$
\hbox{we let}
$$
H(\omega)=\hat \Pi_{12}\hat\Pi_{34}-\hat\Pi_{13}\hat\Pi_{24}+\hat\Pi_{14}\hat\Pi_{23}\in R.
$$
\end{defn}

\begin{lemma}
For $0\le q\le n$, let $V'$ be a $q$-dimensional subspace of $\C^n$.
Then the $Q$-linear map $X$ restricted to $V'\otimes_\C Q$ has the
maximal rank $r:=\min\{\k+2,q\}$, i.e., one of the $r\times r$-minors of
matrix $X$ is a nonzero element of $R$.
\label{lem:wronskian}
\end{lemma}

\begin{proof}
For any matrix $M$, denote by $M_{i_1\dots i_a}^{j_1\dots j_b}$ the
$a\times b$ matrix obtained by taking the rows $i_1$, \dots, $i_a$,
columns $j_1$, \dots, $j_b$ of $M$.
Let $U=(u_{ij})$ be an $n\times q$ matrix whose columns form a basis
of $V'$. It suffices to show that % the top left $r\times r$ minor
$D:=\det\bigl((XU)_{1\dots r}^{1\dots r}\bigr)$
% of the $(\k+2)\times q$ matrix $XU$
is not the zero polynomial.
Using the Cauchy-Binet formula we have
\begin{align*}
D&=\sum_{n\ge h_1>\cdots>h_r\ge1}
\det X_{1\dots r}^{h_r\dots h_1}
\det U_{h_r\dots h_1}^{1\dots r} % \\ &
=% (-1)^{r(r-1)/2}\biggl(\prod_{1\le i<j\le r}(x_i-x_j)\biggr)
\Delta'
\sum_{n-r\ge\lambda_1\ge\cdots\ge\lambda_r\ge0}
c_\lambda s_\lambda(x_1,\dots,x_r),
\end{align*}
where $\Delta'$ is $(-1)^{r(r-1)/2}$ times the Vandermonde determinant
$\prod_{1\le i<j\le r}(x_i-x_j)$,
$s_\lambda=s_\lambda(x_1,\dots,x_r)$ is the Schur function for the partition
$\lambda=(\lambda_1,\lambda_2,\dots,\lambda_r)$, \cite{Macdonald}, and
$c_\lambda=\det U_{h_r\dots h_1}^{1\dots r}$
with $\lambda_i=h_i-(r-i+1)$.
Since the $s_\lambda$ are linearly independent over $\C$,
and since some $c_\lambda\in\C$ do not vanish, $D$ is not the zero polynomial.
\end{proof}

\begin{theorem}\label{thm:params}
The following four conditions for $\omega\in\bigwedge^\k\C^n$ are equivalent:
\begin{enumerate}
\item $\omega$ is decomposable;
\item $\bwedge^\k X(\omega)\in\bigwedge^\k R^{\k+2}$ is decomposable;
\item $G(\omega)\in\bigwedge^2R^4$ is decomposable;
\item \label{Hcondition}$H(\omega)=0$.
\end{enumerate} 
\end{theorem}

\begin{proof}
Since $H$ is the pull-back by $G$ of the Pl\"ucker relation for $\Gamma^2\C^4$,
conditions 3 and 4 are equivalent.  Since $G$ is a GCP map, 1 implies 3.
So we have only to prove that 3 implies 2 and 2 implies 1.

We first prove that 3 implies 2.
Suppose $\alpha:=\bwedge^\k X(\omega)$ is indecomposable.
Then $\beta:=\delta(\alpha)\in\bigwedge^2\C^{\k+2}$ is also indecomposable
since $\delta(\beta)=\alpha$ and $\delta$ is a GCP map.
Hence for some $P_{A,B}\in\P(2,\k+2)$ we have
\begin{equation}
P_{A,B}(\beta)=
\Pi_{i_0i_1}\Pi_{i_2i_3}-\Pi_{i_0i_2}\Pi_{i_1i_3}+\Pi_{i_0i_3}\Pi_{i_1i_2}
\neq0,
\label{foo7}
\end{equation}
where $A=\{i_0\}$ and $B=\{i_1,i_2,i_3\}$ with $1\le i_1<\cdots<i_4\le\k+2$,
and $\beta=\sum_{1\le i<j\le \k+2}\Pi_{ij}e_i\wedge e_j$.
Permutation of the parameters
$(x_i)\mapsto(x_{\sigma(i)})$, where $\sigma\in\mathfrak S(p+2)$, does not
affect the indecomposability of $\alpha$ and $\beta$, and it yields
the action $\Pi_{ij}\mapsto\ve(\sigma,i,j)\Pi_{\sigma(i)\sigma(j)}$
on the coordinates $\Pi_{ij}$ of $\beta$.  Here $\ve(\sigma,i,j)=\pm1$
depends on $i$ and $j$, but in such a way that its net effect on each term
of \eqref{foo7} is the same for all three terms.
Hence for all choices of
$(i_0,\dots,i_3)$, $P_{A,B}(\beta)$ are the same polynomial
upon renaming the variables $x_i$ and a possible change in sign.
In particular, we can take $(i_0,i_1,i_2,i_3)=(1,2,3,4)$, i.e.,
$P_{\{1\}\{2,3,4\}}(\beta)=P_{\{1\}\{2,3,4\}}(\bwedge^2 Z(\beta))\neq0$.
Thus $\bwedge^2 Z(\beta)=G(\omega)$ is indecomposable, proving 3 implies 2.

Now we prove that 2 implies 1.  Suppose $\omega$ is indecomposable.
We will prove the indecomposability of $\bwedge^\k X(\omega)$.

First consider the case where $\omega$ is the sum of \textit{two} decomposable
elements, $\omega =v_1\wedge\cdots\wedge v_\k+w_1\wedge\cdots\wedge w_\k$.
Setting
$V=\langle v_1,\ldots,v_\k\rangle_Q$ and $W=\langle w_1,\ldots,w_\k\rangle_Q$
we have, as seen in the proof of Theorem~\ref{thm:existsABC}, that
$\dim V=\dim W=\k$ and $\dim V\cap W\le \k-2$, so that $\dim V+W\ge \k+2$.
Using Lemma~\ref{lem:wronskian} with $V'=V$,~$W$ and $V+W$ respectively,
we have
$\dim XV=\dim XW=\k$ and $\dim(XV+XW)=\dim(X(V+W))=\k+2$, so that
$\dim(XV\cap XW)=\k-2$.  Hence
$\bwedge^\k X(\omega)= Xv_1\wedge\cdots\wedge Xv_\k + Xw_1\wedge\cdots\wedge Xw_\k$
is indecomposable.

Next we study the general case by induction on $n$. If $n=4$, then $\det X$ is
a nonzero element of $R$, so that 2 implies 1 is obvious.
Suppose the assertion holds for $n-1$.
Let $\omega=\omega_1\wedge e_n+\omega_2$ with
$\omega_1\in\bigwedge^{\k-1}\C^{n-1}$ and $\omega_2\in\bigwedge^\k\C^{n-1}$,
where we regard $\C^{n-1}\subset\C^n$ as in the proof of
Theorem~\ref{thm:existsABC}.
If both $\omega_1$ and $\omega_2$ are decomposable, it reduces to the case
studied above.

Suppose $\omega_1$ is indecomposable. 
Let $X_0$ be the $(\k+1)\times(n-1)$ matrix obtained from $X$ by removing the
last row and column.  Since $Xe_n=\sum_{i=1}^{\k+2}x_i^{n-1}e_i$, we have
\begin{equation}
\begin{split}
\bwedge^\k X(\omega)&=\bwedge^{\k-1}X(\omega_1)\wedge Xe_n+\bwedge^\k X(\omega_2)\\&
=(-1)^{\k+n}x_{\k+2}^{n-1}\bwedge^{\k-1}X_0(\omega_1)\wedge e_{\k+2}+
\text{lower degree terms in $x_{\k+2}$}.
\end{split}
\label{Xpandinx_\k+2}
\end{equation}
By the induction hypothesis $\bwedge^{\k-1}X_0(\omega_1)$ is indecomposable,
so the right-hand side of \eqref{Xpandinx_\k+2}
is also indecomposable, as seen by expanding the
Pl\"ucker relations $P_{A,B}(\bwedge^\k X(\omega))$ with $A\cap B\ni \k+2$
in powers of $x_{\k+2}$ and taking the coefficients of $x_{\k+2}^{2n-2}$.

Suppose $\omega_1$ is decomposable and $\omega_2$ is not. Then
$\bwedge^{\k-1} X(\omega_1)$ is decomposable, and by the induction hypothesis
$\bwedge^\k X(\omega_2)$ is not.
We prove by contradiction that $\bwedge^\k X(\omega)$ cannot be decomposable:

Let $\p'_{i_1,\dots,i_{\k-1}}$ (resp.\ $\q_{i_1,\dots,i_\k}$) be the
coordinates of $\bwedge^{\k-1} X(\omega_1)$ (resp.\ $\bwedge^\k X(\omega_2)$). Thus
$\bwedge^\k X(\omega_1\wedge e_n)=\bwedge^{\k-1} X(\omega_1)\wedge Xe_n$ has the
coordinates
\begin{equation}
\p_{i_1,\dots,i_\k}=
\sum_{\nu=1}^\k(-1)^{\k-\nu}\p'_{i_1,\dots,i_\k\backslash i_\nu}x_{i_\nu}^{n-1}.
\label{PcoordExpn}
\end{equation}
By assumption, $\{\p_{i_1,\dots,i_\k}\}$ satisfies all of the Pl\"ucker relations
in $\P(\k,\k+2)$, and $\{\q_{i_1,\dots,i_\k}\}$ does not satisfy some relation,
say $P_{A,B}$, in $\P(\k,\k+2)$.  Here,
using the symmetry argument as we used in the proof of ``3 implies 2'' above,
we may assume $A=\{1,5,6,\dots,\k+2\}$ and $B=\{2,3,4,5,6,\dots,\k+2\}$.
Denoting the sequence of indices $5,6,\dots,\k+2$ by $\lambda$, and hence
$\p_{ij\lambda}=\p_{ij56\dots \k+2}$ etc., we have thus
\begin{align}
\p_{12\lambda}\p_{34\lambda}-\p_{13\lambda}\p_{24\lambda}
+\p_{14\lambda}\p_{23\lambda}&=0,\label{PluckPP}\\
\q_{12\lambda}\q_{34\lambda}-\q_{13\lambda}\q_{24\lambda}
+\q_{14\lambda}\q_{23\lambda}&\neq0.\label{PluckQQ}
\end{align}
If $\bwedge^\k X(\omega)$ is decomposable, then
$\{\p_{i_1,\dots,i_\k}+\q_{i_1,\dots,i_\k}\}$, the coordinates of
$\omega$, must satisfy all the relations in $\P(\k,\k+2)$ and so the above
$P_{A,B}$ in particular; thus by using \eqref{PluckPP} we have
\begin{equation}
\begin{split}
0&=(\p_{12\lambda}+\q_{12\lambda})(\p_{34\lambda}+\q_{34\lambda})-(\p_{13\lambda}+\q_{13\lambda})(\p_{24\lambda}+\q_{24\lambda})\\
&\hphantom{{}={}}+(\p_{14\lambda}+\q_{14\lambda})(\p_{23\lambda}+\q_{23\lambda})\\
&=(\p_{12\lambda}\q_{34\lambda}-\p_{13\lambda}\q_{24\lambda}+\p_{14\lambda}\q_{23\lambda})
+(\q_{12\lambda}\p_{34\lambda}-\q_{13\lambda}\p_{24\lambda}+\q_{14\lambda}\p_{23\lambda})\\
&\hphantom{{}={}}+(\q_{12\lambda}\q_{34\lambda}-\q_{13\lambda}\q_{24\lambda}+\q_{14\lambda}\q_{23\lambda}).
\end{split}
\label{PluckP+Q}
\end{equation}
By the definition of linear map $X$, $\q_{ij\lambda}$ are polynomials in $x_r$,
$r=i$,~$j$, $5,\dots,\k+2$, with no $x_r$ having the $(n-1)$st or higher power,
and $\p_{ij\lambda}$ are polynomials in the same set of $x_r$ with $x_r^{n-1}$
appearing only as the last factor in each term of \eqref{PcoordExpn}.
Moreover, for each $r\in\{1,\dots,4\}$, each term on the right-hand side of
\eqref{PluckP+Q} is the product of a polynomial which depends on $x_r$ and
another which does not.  Thus comparing the coefficients of $x_r^{n-1}$
($r=2$, 3, 4) on both sides of \eqref{PluckP+Q}, we have
\begin{gather*}
\p'_{1\lambda}\q_{r'r''\lambda}-\p'_{r'\lambda}\q_{1r''\lambda}+\p'_{r''\lambda}\q_{1r'\lambda}=0,\qquad
(r',r'')=(3,4),\,(4,2),\,(2,3).
\end{gather*}
After using these to eliminate $\q_{23\lambda}$, $\q_{24\lambda}$ and $\q_{34\lambda}$,
the left-hand side of \eqref{PluckQQ} becomes
$$
\q_{12\lambda}\frac{\p'_{3\lambda}\q_{14}-\p'_{4\lambda}\q_{13}}{\p'_{1\lambda}}
-\q_{13\lambda}\frac{\p'_{2\lambda}\q_{14\lambda}-\p'_{4\lambda}\q_{12\lambda}}{\p'_{1\lambda}}
+\q_{14\lambda}\frac{\p'_{2\lambda}\q_{13\lambda}-\p'_{3\lambda}\q_{12\lambda}}{\p'_{1\lambda}}=0,
$$
which is a contradiction. Hence $\bwedge^\k X(\omega)$ is indecomposable,
completing the proof of ``2 implies 1''. \end{proof}

\section{Concluding Remarks}

We obtained a set $\P'(\k,n)$ of rank 6 quadratic forms on $\bigwedge^\k\C^n$,
which is in general much smaller than the set of standard Pl\"ucker relations,
and yet is capable of determining the decomposable elements in
$\bigwedge^\k\C^n$.
Every element of $\P'(\k,n)$ is obtained by pulling back the one nontrivial
Pl\"ucker relation on $\bigwedge^2\C^4$.  This means:

1.~The elements of $\P'(\k,n)$ define
quadric hypersurfaces which are isomorphic to each other by a GCP map,
and the Grassmannian is obtained as the intersection of those isomorphic
quadrics.  This may be of interest to geometers who have already
noted and used the previously known fact that this was true when $\k=2$
\cite{rk6ref1,rk6ref2}.

2.~The 3-term Pl\"ucker relation is in a sense ``universal'', providing an explanation for the special role played by the 3-term relation
in applications like soliton theory \cite{Duzhin,GK}.  In this respect,
condition~\eqref{Hcondition} in Theorem~\ref{thm:params} should be viewed as
being analogous to the use of a parameter-dependent Fay-Hirota type difference
equation to characterize KP tau functions \cite{Miwa,Tak-Tak}.

3.~The Grassmannian $Gr(\k,n)$ is the  
intersection of the pullbacks of $Gr(2,4)$ under all GCP maps from  
$\bigwedge^\k\C^n$ to $\bigwedge^2\C^4$.

% Geometrically, the results presented above also prove that the Grassmannian variety $Gr(\k,n)$ can be constructed as the intersection of a collection of isomorphic quadric hypersurfaces.  In particular, we note that all the relations in $\P'(\k,n)$ have isomorphic vanishing loci as they are mapped to each other by the induced map of a linear map on $\C^n$.  This may be of interest to geometers who have already noted and used the previously known fact that this was true when $\k=2$ \cite{rk6ref1,rk6ref2}.

%  Finally, we note that our main results help to clarify the method which is being implicitly used  in soliton theory when the Hirota Bilinear Difference Equation (HBDE) characterizes KP tau-functions  \cite{GK,Miwa}. 
%  In particular, condition \eqref{Hcondition} in Theorem~\ref{thm:params} should be viewed as being analogous  to the use of one functional equation (the HBDE)  to replace the many equations of the KP hierarchy corresponding directly to individual Pl\"ucker relations \cite{KP-Plucker1,KP-Plucker2}.

% \par\bigskip\par\noindent\textbf{Acknowledgments:}
\subsection*{Acknowledgments:}
KP and AR worked on this paper as part of a summer undergraduate research experience at the College of Charleston, with KP receiving financial support from the Math Department and AR receiving financial support from the ``SURF''.  We are grateful for assistance and advice from 
Malcolm Adams, 
Michael Gekhtman,  
Trygve Johnsen,
Mitch Rothstein, 
Oleg Smirnov,
Robert Varley and the referee.
% The proof of Lemma~\ref{lem:subspaces} in earlier version of the paper
% had an error, which was then found by the referee.

\end{document}